\documentclass[12pt]{article}
\usepackage{amsmath,amsthm,amssymb}

\theoremstyle{plain}
\newtheorem{theorem}{Theorem}
\newtheorem{lemma}{Lemma}
\newtheorem{definition}{Definition}

\newtheorem{definition-theorem}{Definition-Theorem}

\theoremstyle{definition}
\newtheorem{remark}{Remark}

\title{Intrinsic metric on graded graphs, standardness, and
 invariant measures}
\author{A.~M.~Vershik\thanks{%
St.~Petersburg Department
 of Steklov Mathematical Institute, e-mail: \mbox{vershik@pdmi.ras.ru}.
 Supported by the RFBR grant 13-01-12422-ofi-m. The author also thanks the Hausdorff Institute for Mathematics in Bonn
for support during the semester ``Universality and Homogeneity''
(September--December 2013).}
}

\begin{document}

  \maketitle

\begin{abstract}
We define a general notion of a smooth invariant (central) ergodic measure on the space of paths of
an $N$-graded graph (Bratteli diagram).  It is based on the notion of standardness of the
tail filtration in the space of paths, and the smoothness criterion uses the so-called
intrinsic metric which can be canonically defined on the set of vertices of these
graphs. In many cases known to the author, like the Pascal graph, the Young graph, the space of configurations,
all ergodic central measures are smooth (in this case, we say that the graph is smooth).
But even in these cases, the intrinsic metric is far from being obvious and does not coincide with the ``natural''
metric. We apply and generalize the theory of filtrations developed by the author during the last forty years to the case
of tail filtrations and, in particular,  introduce the notion of a standard filtration as a generalization
to the case of semi-homogeneous filtrations of the notion of a standard homogeneous (dyadic) filtration in the sense
of that theory. The crucial role is played by the new notion of \textbf{intrinsic
semi-metric on the set of vertices of a graph} and the notion of regular paths, which allows us to refine
the ergodic method for the case of  smooth measures. In future,
  we will apply this new approach to the theory of invariant measures in combinatorics, ergodic theory, and the theory of C$^*$-algebras.
\end{abstract}


\section{The simplex of invariant measures and the ergodic method}

Let $\Gamma$ be an $N$-graded graph, or, more exactly, a ``Bratteli diagram\footnote{This means that there is one vertex $0$ of degree $0$, finitely many
vertices of each given degree $n$, and each vertex at levels $n>0$ has at least one
predecessor and at least one follower; beside that, an edge may join only vertices of adjacent levels.
In general, we allow finitely many edges between two vertices, but if this is not mentioned explicitly, we consider
edges without multiplicities. We do not use the term ``Bratteli diagram,'' because we call the vertices of $\Gamma$
``diagrams,'' and paths ``tableaux,'' just as in the theory of the Young graph.}.''
The set of vertices of level $n$ is denoted by $\Gamma_n$.
A finite path is a sequence of edges $(e_1,e_2,\dots, e_n)$ such that the end vertex of $e_k$ is the initial vertex of  $e_{k+1}$, $k=1, \dots, n-1$. In the case where the graph has no multiple edges, we may say that a path is a sequence of diagrams (vertices).
A finite path that starts at $0$ and ends at a vertex $v \in \Gamma$ will be called  a ``tableau with diagram $v$.''
A maximal path,  or an infinite tableau, $t$ is an infinite sequence of
vertices $t_n$, $n=0,1, \dots$, where $t_0=\varnothing$. Again, if the graph
$\Gamma$ has no multiple edges, then a path $\{t_n\}$ is a sequence of vertices (diagrams).

Let  $T(\Gamma)$ be the space of all maximal infinite tableaux. There is a
natural structure of inverse limit, or topological Cantor space, on the space $T(\Gamma)$, defined via
clopen,  or cylinder, sets; we denote by ${\mathfrak  A}_n$ the $\sigma$-field of cylinder sets of level
$n$. Thus we have a decreasing sequence of $\sigma$-fields $\{{\mathfrak  A}_n\}$, $n \in \mathbb  N$,
 which we call the {\it tail filtration} on the space $T(\Gamma)$. The language of $\sigma$-fields is parallel to that of measurable partitions: we will consider decreasing sequences of partitions $\{\xi_n\}$ corresponding to  $\{{\mathfrak  A}_n\}$, where
 $\xi_n$ is the partition of the space $T(\Gamma)$ into classes of paths that coincide after level $n$;
 each block (element) of  $\xi_n$ has finitely many points. The set-theoretical intersection
 of the partitions $\{\xi_n\}$ is the {\it tail partition}, which defines the {\it tail equivalence relation}: the classes of the tail equivalence relation, which we denote by $\xi_{\infty}$, are unions of increasing classes of the partitions $\xi_n$ over $n$.
  (We use the same symbol $\xi_{\infty}$ also for the limit partition.)

   In most interesting cases, $\xi_{\infty}$ is ``ergodic,'' which means that there is no countable set of Borel functions that are constant on all classes of the relation and can distinguish these classes. Each such class $C \in \xi_{\infty}$ is a countable set, and we can restrict the sequence of partitions $\{\xi_n\}$ to this element $C$. This is what we called the ``hierarchy,''  or the filtration on the element $C$; denote it by $\{\xi^C_n\}$.  If we want to emphasize that an element of $\xi^C_n$ is the set of all tableaux with diagram $v\in \Gamma_n$, we denote it by $\xi_n(v)$.
   We refer the reader to the theory of filtrations (or decreasing sequences of measurable partitions; see, e.g., \cite{V4}) and the
   notion of standard filtrations. This notion was defined for homogeneous filtrations, but the {\it tail partition} for
   the space of paths with a central measure is not homogeneous, but  what we call semi-homogeneous, and we need to
   give a generalization of the notion of standardness; this will be done in the paper mentioned below,
   but we will not use this notion here.

  The ``dimension,'' or ``binomial coefficient,'' $\dim u$ of a vertex (diagram) $u\in\Gamma_n$ is the number of paths that start at $\varnothing$ and end at $u$.\footnote{Note that in the case of multiple edges, a path is a sequence of edges rather than vertices.} The {\it probability measure $\nu_v$ induced by a diagram $w$ at the previous level} is defined as follows: $$\nu_v(w)=\frac{\dim w}{\dim v},$$ where $w$ runs over all diagrams of level $n-1$ that precede
  $v$ (which we denote as $v\succ w$); recall that $\dim v=\sum_{w\prec v}\dim w$.\footnote{Instead of this measure, we could consider other averages of $\dim v$, e.g., the  root mean square etc., which are similar to the Plancherel measure; however, it is not clear whether this is useful.}

\def\sim{\operatorname{sim}}
\def\Inv{\operatorname{Inv}}
\def\Ex{\operatorname{Ex}}

   \begin{definition}[Projection and the inverse limit]
   Denote the simplex of all formal convex combinations of vertices of $\Gamma_n$ by $\sim(\Gamma_n)$ and define an
   affine projection $p_n: \sim(\Gamma_{n+1}) \rightarrow  \sim(\Gamma_n)$ as the convex extension of the
   projection of vertices: $p(v)= \sum_{w:w\prec v} \nu_v(w) \delta_w$, where $v \in \Gamma_{n+1}, w\in \Gamma_n$.

   Denote the corresponding inverse limit by  $$\Inv(\Gamma)\equiv \lim_{\leftarrow}\{\sim(\Gamma_n), p_n\}.$$
   \end{definition}

We also consider the projections $p_{n,m}=p_n p_{n-1,m}\dots p_{m+1,m}$ and the corresponding limit $p_{\infty,m}:
   \Inv(\Gamma)\rightarrow \sim(\Gamma_n)$.

   \begin{remark} {\rm The projections $p_{n,m}=p_np_{n-1,m}\dots p_{m+1,m}$ and $p_{\infty,m}$ may be not epimorphisms onto $\mbox{sim}(\Gamma_m)$ for all $m$. We can represent the same limit  $\mbox{Inv}(\Gamma) $ as the inverse limit of the epimorphic images, which we denote by $\gamma_n\equiv p_{\infty,n}\Inv(\Gamma) \subset\sim(\Gamma_n)$:  $$ \lim_{\leftarrow}(\gamma_n, p_n)\equiv \Inv(\Gamma).$$
     }\end{remark}

   \begin{lemma}
   The inverse limit $\Inv(\Gamma)$ is an affine simplex (with the inverse limit topology) that is canonically isomorphic to the simplex of all central invariant probability measures on the space $T(\Gamma)$ of paths on the graph $\Gamma$ with the weak topology on the space of measures on $T(\Gamma)$.
   The extreme points $\mu \in \Ex \Inv(\Gamma)$ are the ergodic central measures, and each of them is the limit of a sequence of points $\nu_{n_k} \in \Inv(\Gamma)$ whose projection to $\sim(\Gamma_{n_k})$ is an extreme point of the image, $k=1,2, \dots $.
   \end{lemma}

    The proof of this lemma is based on the usual standard arguments of convex affine geometry. Thus our problem is to find
    the set of extreme points of an inverse limit of simplices. There are three kinds of difficulties.

     1) Note that the inverse limit $$\lim_{\leftarrow}\{\sim(\Gamma_n), p_n\}$$
      of finite-dimensional simplices does not allow one to represent any ergodic central measures from $\Inv(\Gamma)$
     as a limit of extreme points of finite simplices; this is an ``approximation from outside,'' because we have no
     embedding of these points to $\Inv(\Gamma)$.

     2) In order to find an approximation of ergodic measures (= extreme points), we usually apply the ``ergodic method''
     (see below). But if we use generic points to find the parameters of the limit measures (the joint distributions, correlation functions, and so on),  the calculations can be very difficult and do not allow one to discover important properties of the measure $\mu$.

     3) The set of ergodic measures $\Ex \Inv(\Gamma)$ can be  dense in $\Inv(\Gamma)$ -- this is a so-called Poulsen
      simplex; in this case, there is no good parametrization of the ergodic measures in principle, and the problem is not well
      posed. We must have a criterion how to separate this case from the smooth case in which a description is possible: the inverse limit
      is a Bauer simplex, in which, by definition, the set of ergodic measures is closed.

    The application of ergodic ideas can be briefly described as follows. In order to find all ergodic measures, we can use the ergodic theorem; or, more exactly, the martingale theorem asserts that if $\mu$
   is an arbitrary ergodic measure on $T(\Gamma)$, then the values of  $\mu$ on cylinder sets $C$ can be calculated as follows. There is a set $T_0(\Gamma)\subset T(\Gamma)$ of full $\mu$-measure such that for each path $t\in T_0(\Gamma)$ there exists a limit $$\mu(C)=\lim_n\frac{|\{s=(s_0=\varnothing, s_1, \dots, s_n=t_n)\in C\}|}{N(n)},$$
   where $s=(s_0,s_1,\dots, s_n=t_n)$ is an arbitrary finite path with the last vertex $t_n=s_n$ and $N(n)=\dim t_n$
   (i.e., $\mu(C)$ is the limit of the fraction of paths that end at $t_n$ and belong to $C$). Thus it is possible in principle to
   describe all invariant measures using the ergodic theorem in this way.
   This method of describing the invariant measures was called the ``ergodic method,'' and it is indeed  very powerful, see
   \cite{V1}, and also \cite{V2,V3} for many further examples of the application of this idea.

    But in concrete situations, we must use more specific properties of spaces and measures, and refine the choice of
    ``generic points.'' In the case of theorems similar to Aldous' theorem \cite{A}, we use other properties (see \cite{V2,V3}).

    In this paper, we suggest a new elaboration of the ergodic method for the case of central measures; in this situation, there is a more specific choice of paths for which the calculations reduce to a pure combinatorial problem. The main idea came
    from the theory of filtrations; more exactly, it is a combinatorial version of the idea of towers of measures \cite{V4}, which in turn used the ``Kantorovich functor'' \cite{K,V5}. We hope that this will help to obtain new theorems on the description of the list of invariant measures, and to give new proofs for the known cases in which the existing proofs  are rather cumbersome, like the proof of Thoma's theorem for the Young graph.

     \section{The intrinsic metric on the graph}

Assume that we have a finite metric space $(X,d)$ with metric $d$. We define the transportation (Kantorovich) metric on the simplex of all probability measures on $X$  as follows. Let $\nu_1=\sum_{i\in I} a_i\delta_{x_i}$, $\nu_2= \sum_{j\in J}b_j\delta_{y_j}$ be two measures; here  $x_i,y_j \in X$, $\sum a_i=\sum b_j=1$, $a_i\geq 0$, $b_j \geq 0$. Then
  $$k_d (\nu_1,\nu_2)=\min_{\{C=\{c_{i,j}\}\}}\sum_{i,j} d(x_i,y_j)c_{i,j},$$
   where
   $$ C=\{\{c_{i,j}\}:\sum_{i\in I} c_{i,j}=b_j, j\in J;  \sum_{j\in J} c_{i,j}=a_i, i\in I; c_{i,j}\geq 0\}.$$

\def\Vert{\operatorname{Vert}}

Below we will consider a graph without multiple edges, so a path is a sequence of vertices. Let $\Vert(\Gamma)=\bigcup_{n=0}^{\infty} \Gamma_n$
be the  set of all vertices of~$\Gamma$. Now we want to define a canonical ``intrinsic semimetric'' on the set
  $\Vert(\Gamma)$ of all vertices of an arbitrary graph $\Gamma$ -- a Bratteli diagram.

\begin{definition}[Intrinsic semimetric on the set $\Vert(\Gamma)$ of all vertices of a graded graph $\Gamma$]
First we define by induction the ``intrinsic semimetric'' on the set of all vertices of a given level
$\Gamma_n$, $n=1,2, \dots $.

 There is one vertex of level $0$, so the space $\Gamma_0=\{\varnothing\}$ is the one-element metric space.

 Define a metric $\rho$ on the finite space of vertices of the first floor $\Gamma_1$ by the formula $\rho (u,v)=1$ if $u\ne v$.

Assume that we have already defined the intrinsic semimetric $\rho$ on $\Gamma_n$, $n\geq 1$. Let $u,v$ be two vertices (diagrams) at the level $n+1$.
   Then the semimetric $\rho$ on this level is defined as follows: $$\rho (u,v)=k_{\rho_n}(\nu_u,\nu_v).$$

 If $u\in \Gamma_n$ and $w\in \Gamma_{n+1}$, then $$\rho(u,w)=k_{\rho_n}(\nu_u, \delta_w).$$

 Finally, let  $v \in\Gamma_n$, $z\in \Gamma_m$ with $m<n$; then $$\rho(v,z)=\min_{\{w_k\}} \sum_{k=m}^{n-1} \rho(w_k,w_{k+1}),$$
 where a sequence $\{w_m=z, w_2,\dots, w_n=v\}$ is a path from $z$ to $v$, and the minimum is taken over all such paths.

 Thus we have defined the {\it intrinsic semimetric $\rho=\rho_{\Gamma}$} on $\Vert(\Gamma)$.
 \end{definition}

Note that the distance between two diagrams of the same level whose sets of preceding vertices coincide ($\nu_u=\nu_v$), or between vertices $u,w$ such that $w$ is the unique predecessor of $u$, vanishes, so that in general $\rho$ is a semimetric, not a true metric. But in many cases, $\rho$ is indeed a metric.

 The initial metric $\rho_1$ on $\Gamma_1$ can be replaced by some other metric; moreover, we can also start
 from a metric on a level $\Gamma_n$ as the initial metric. The asymptotic properties of the intrinsic semimetric, as we will see, do not depend on these changes.

\begin{remark}{\rm	
 The semimetric space $(\Vert(\Gamma), \rho)$ is not an inductive limit of the spaces $(\Gamma_n, \rho_n)$, because there are no embedding isometries $\Gamma_n \to\Vert(\Gamma)$.
 }\end{remark}

\begin{definition}
Let $\widehat \Vert(\Gamma)$ be the quotient of the completion  of the space $(\Vert(\Gamma), \rho)$ by the equivalence relation into the classes of points with zero distance. We  call it the  \textbf{``intrinsic limit space''} with \textbf{``intrinsic metric''} $\hat \rho$.
\end{definition}

The {\it accumulation set} $I(\Gamma)\equiv \widehat \Vert(\Gamma) \setminus \Vert(\Gamma)$ is the goal of our considerations. There is no natural way to consider vertices of $\Gamma$ as points of $I(\Gamma)$ (see Remark 2 above). But we may say that a point of $I(\Gamma)$ is the limit, in the sense of the intrinsic metric, of some path of  $T(\Gamma)$ regarded as a sequence of vertices.

A path (regarded as a sequence of vertices) is called {\it regular} if the set of its vertices
  is a Cauchy sequence in $(\Vert(\Gamma),\rho)$ with respect to the intrinsic metric $\rho$, or, equivalently, a convergent sequence in $I(\Gamma)$. Denote by $T_{{\rm reg}}(\Gamma)$ the subset of regular paths in the compact space  $T(\Gamma)$ of all paths. Each point of $I(\Gamma)$ is the limit of a regular path.

   Our main result is the following theorem.

\begin{theorem}
   Every regular path canonically determines an ergodic central measure on the space $T(\Gamma)$; this measure
   depends only on the limit of the path regarded as a point of $I(\Gamma)$; thus we have a map $$I(\Gamma) \rightarrow \Inv(\Gamma),$$
   which is a monomorphism. In the case of a smooth graph or a standard tail filtration (see below), this map is an epimorphism, so we have a method of describing the list of all ergodic invariant measures on $T(\Gamma)$ for smooth graphs.
\end{theorem}

  This theorem provides a refinement of the ergodic method, which consists in a very precise selection of paths that generate invariant measures.

  The proof of the first part, the existence of a central measure generated by a regular path, is more or less obvious, by the definition of the intrinsic metric and the fact that the path is regular. The ergodicity of the limit measure also follows from the convergence of the path; a detailed proof will be published elsewhere.

The most important thing is that under the limit measure, as a central measure on the space of paths $T(\Gamma)$, the tail filtration on $T(\Gamma)$ is standard (see the corresponding remarks above). The standardness implies a stronger (compared with the usual one) form of the theorem
on convergence of reverse martingales. This  strengthening is closely related to
the substitutional ergodic thorems in the sense of
 \cite{V6}.

   We can conclude that the standardness leads to the smoothness of the graph.

\section{Comments and examples}

\subsection{Some comments}

{\bf1.} The definition of the intrinsic metric can be given in
terms of linear algebra only. Indeed, we define a metric on the set of
linear subspaces of a linear space, but these subspaces have an
additional structure, namely, the decomposition into a sum of
subsubspaces, each of which has a further decomposition, etc.; these
decompositions end with one-dimensional subspaces  whose unit vectors constitute the
Gelfand--Tsetlin basis. The definition of the intrinsic metric is based on a {\it
correspondence between subspaces} that takes into account the
hierarchical structure.

\smallskip\noindent
{\bf 2.} The notion of a semi-homogeneous standard filtration is based on the {\it standardness criterion} for homogeneous
(for example, dyadic) filtrations, namely, it uses the iteration of the Kantorovich metric. In the definition above, we also
implicitly use the same idea. That is why the compactness of $\Vert(\Gamma)$ with respect to the intrinsic metric
is equivalent to the standardness.

\smallskip\noindent
{\bf 3.} One can define many various similar notions of metric even for a given diagram: we can use noncentral measures
instead of central ones, a cocycle on the tail equivalence relation, etc.

\smallskip\noindent
{\bf4.} A Bratteli diagram generates a Markov compactum (in general, nonstationary), so we can pose many new questions
about Markov processes and random walks on graphs and groups. The set of central measures in this case is the
 {\it exit boundary},  or the Poisson--Furstenberg boundary. We hope that our considerations will provide a new method of
 calculating it.

\smallskip\noindent
{\bf 5.} The intrinsic metric explains the difference between ``Poulsen'' and ``Bauer'' simplices (or the smooth and nonsmooth cases) which we mentioned at the beginning of the paper.
 See also the paper \cite{7}, where we consider the general problem of
description of traces on $C^*$-algebras.

\subsection{Some examples}

{\bf 1.}  For the usual Pascal triangle, the set $\Gamma_n$  of vertices of level $n$, $n=1,2, \dots$, equipped with the intrinsic metric is isometric
to the set of $n+1$ points $\frac {i}{n}$, $i=0,1, \dots, n$, of the interval with  the usual metric.
Therefore, the natural limit of $\Gamma_n$ regarded as a metric space is the unit interval with the usual metric.

Recall that the space of all maximal paths of $\Gamma$ is the Cantor space $\{0;1\}^{\infty}$, which can be mapped
to $[0,1]$, but this isomorphism with the space of paths $T(\Gamma)$ (except a countable subset) has nothing to do with our limit: if we consider a path as a sequence of points of $\Gamma_n$, then we restrict ourselves to a small part of paths that have limits.

By well-known theorems (for instance, de Finetti's theorem),  the set of ergodic invariant measures on $\Gamma$ is also parameterized by the interval $[0,1]$ (if we distinguish the states $0$ and $1$). Thus our limit can be interpreted as a mysterious correspondence between regular paths and ergodic measures (in the sense that each regular path generates an ergodic measure by the LLN or ergodic theorem).

\smallskip\noindent
{\bf 2.} For the Pascal tetrahedron (in dimension 3), calculations show that the set of vertices of level $n$ equipped with the intrinsic metric is isometric to the set of all points of the plane simplex with coordinates of the form $\frac{i}{n}$ equipped with the
hexagonal metric, generated by the Banach norm whose unit ball centered at the barycenter of the simplex is a regular hexagon. The same corollary as in the first example about an isomorphism of the set of paths and a real simplex can also be derived in this case. A similar answer is true for the Pascal graph of arbitrary dimension.

\smallskip\noindent
{\bf3.} For the Young graph, the limit of $\Gamma_n$ with the intrinsic metric is supposedly the Thoma simplex,
and this should give a combinatorial proof of the theorem on central measures on the Young graph.
The main lemma asserts that the convergence in the intrinsic metric coincides with the convergence of the ``{\bf frequencies} of rows and columns.''

\end{document}